\documentclass[a4paper,12pt]{amsart}

\usepackage{amssymb,amsmath,amscd}

\def\binom#1#2{\left( \begin{array}{c}{#1} \\ {#2}\end{array} \right) }

\newcommand{\N}{\mathbb{N}}
\newcommand{\R}{\mathbb{R}}

\newcommand{\bk}{\mathrm{\bf k}}  
\newcommand{\br}{\mathrm{\bf r}}
\newcommand{\bs}{\mathrm{\bf s}}
\newcommand{\bl}{\mathrm{\bf l}} 
\newcommand{\bt}{\mathrm{\bf t}}

\newcommand{\ba}{\mathrm{\bf a}}
\newcommand{\bm}{\mathrm{\bf m}}

\newcommand{\Lie}{\mathcal{L}}       

\newtheorem{thm}{Theorem}[section]
\newtheorem{cor}[thm]{Corollary}

\theoremstyle{definition}
\newtheorem{defin}[thm]{Definition}
\newtheorem{rem}[thm]{Remark}

\numberwithin{equation}{section}  

\begin{document}

\title{On the multilinear Hausdorff  problem of moments}

\author{A. Ibort}
\address{Depto. de Matem\'aticas, Univ. Carlos III de
Madrid, Avda. de la Universidad 30, 28911 Legan\'es, Madrid, Spain.}
\email{albertoi@math.uc3m.es}

\author{P. Linares, J.G. Llavona}
\address{Departamento de An{\'a}lisis Matem{\'a}tico, Facultad de Matem{\'a}ticas, 
Universidad Complutense de  Madrid, 28040 Madrid, Spain}
\email{plinares@mat.ucm.es, jl\_llavona@mat.ucm.es}

\keywords{Problem of moments, multilinear, polymeasures, second order stochastic processes.}

\begin{abstract}
Given a multi-index sequence $\mu_{\bk}$, $\bk = (k_1, \ldots, k_n) \in \N_0^n$,
necessary and sufficient conditions are given for the existence of
a regular Borel polymeasure $\gamma$ on the unit interval $I= [0,1]$ such that
$\mu_{\bk} = \int_{I^n} t_1^{k_1}\otimes \cdots \otimes t_n^{k_n} \, \gamma$.   
This problem will be called the weak mul\-ti\-li\-near Hausdorff problem
of moments for $\mu_{\bk}$.   Comparison with
classical results  will allow us to relate the weak multilinear Hausdorff problem with the multivariate Hausdorff problem.
A solution to the strong multilinear Hausdorff problem of moments will be
provided by exhibiting necessary and sufficient conditions for the existence
of a Radon measure $\mu$ on $[0,1]$ such that $L_\mu(f_1,\ldots, f_n) =
\int_{I} f_1(t) \cdots f_n(t) \, \mu (dt)$ where $L_\mu$ is the $n$-linear moment  functional
on the space of continuous functions on the unit interval defined by the sequence $\mu_{\bk}$.
Finally the previous results will be used to provide a characterization of
a class of weakly harmonizable stochastic processes with bimeasures supported on compact sets.
\end{abstract}

\maketitle

MSC Classification:  Primary 44A60; Secondary 46G25.

\section{Introduction}
The moment functional $L_\mu$ associated to a sequence $\mu_k$,
$k\in\N_0$ of real numbers is the element in the (algebraic) dual of
the space of polynomials $\R [t]$ defined by  $L_\mu (p) = \sum_{k \geq 0} p_k \mu_k$,
where $p(t) = \sum_{k \geq 0} p_k t^k \in \R [t]$ is an arbitrary
polynomial.  Given an interval $I \subset \R$ the classical
problem of moments for the sequence $\mu_k$ asks for the
integrality of the linear operator $L_\mu$, that is, under what conditions there exists a (possibly signed) Radon measure $\mu$ on
$I \subset \R$ such that $L_\mu(t^k) =
\int_I t^k d\mu (t)$, $k =0,1\dots$.   If $I$ is the unit interval $[0,1]$ the problem
of moments is known as the Hausdorff moment problem.

The well-known solution to the classical Hausdorff problem (see  for instance
\cite{Sh70}  and references therein) establishes that
such a measure $\mu$ exists provided that
there is a constant $C$ such that:
\begin{equation}\label{linear_bound} \sum_{m = 0}^k | \lambda_{(k;m)} | < C ,
\end{equation}
for all $k = 0,1,\ldots$,  where $ \lambda_{(k;m)} = \binom{k}{m} L_\mu (t^m (1-t)^{k-m}) $.

There are various natural extensions of the moment problem to
the multilinear case.  
Given the multi--index sequence $\mu_{\bk}$,  $\bk = (k_1, \ldots, k_n) \in \mathbb{N}_0^n$, we will denote as before by $L_\mu$ the $n$--linear functional defined on the
space of polynomials $\R [t]$ by:
\begin{equation}\label{functional}
L _\mu (t^{k_1}, \ldots, t^{k_n}) =  \mu_{\bk} \quad \forall \bk = (k_1,\ldots, k_n) \in \mathbb{N}_0^n .
\end{equation}
Thus the strong multilinear Hausdorff problem
of moments will consist in determining under what conditions there exists a (possibly signed) Radon
measure $\mu$ on $[0,1]$ such that
\begin{equation}\label{strong_moment}
L _\mu (p_1(t), \ldots, p_n(t) ) = \int_0^1 p_1(t) \cdots p_n(t)\, \mu (dt) , \quad \forall p_1(t), \ldots, p_n(t) \in \R [t]  .
\end{equation}
A weaker version of the multilinear Hausdorff moment problem, the classical (multivariate) Hausdorff moment
problem, can be stated by demanding the
existence of a (possibly signed) Radon measure $\mu$ on $[0,1]^n \subset \R^n$ such that:
\begin{equation}\label{classical}
L_\mu(p_1(t_1) , \ldots, p_n (t_n) ) = \int_{I^n} p_1(t_1) \cdots p_n (t_n) \, \mu (dt_1, \ldots , dt_n) .
\end{equation}

Integrality properties of bilinear functionals related to the notion of the total variation
of functions in many variables were studied by Morse and Transue \cite{Mo49}. These questions led to the concept of  $\mathbb{C}$--bimeasures \cite{Mo56}.
These results were also deeply rooted in the problem of studying the structure
of stationary stochastic processess
(see for instance \cite{Ra82} and references therein).
The notion of bimeasures, and polymeasures in general, provide a
natural framework to answer these questions.   Thus we are naturally led to
consider an even weaker version of the classical multilinear
Hausdorff problem of moments.   
We will say that $\mu_{\bk}$ satisfies the weak multilinear
Hausdorff problem of moments if there exists a polymeasure
$\gamma$ on $\mathrm{Bo}[0,1] \times \cdots \times
\mathrm{Bo}[0,1]$ such that:
\begin{equation}\label{weak_problem} 
\mu_{\bk} = \int_{I^n} t_1^{k_1} \otimes \cdots \otimes t_n^{k_n} \, \gamma (dt_1, \ldots, dt_n) , \quad \forall \bk = (k_1, \ldots, k_n) \in \mathbb{N}_0^n .
\end{equation}
As it turns out, the solution to the weak multilinear Hausdorff problem of moments
is given by a nontrivial generalization of condition eq. (\ref{linear_bound}) as it will be
proved in section 2, Thm. \ref{cont_multi}.   This condition is different from the characterization obtained in the  analogous 
weak multilinear trigonometric problem of moments \cite{Ib10}.   Using these ideas the classical condition  (\ref{linear_bound}) and the
solution to the classical multivariate Hausdorff problem are easy consequences of the 
general properties of polymeasures on compact sets as it will be discussed briefly at the end of section 2.
The strong multilinear Hausdorff problem of moments will be solved in section 3 by using recent results on integral representations of orthogonally additive polynomials on Banach lattices.  In addition the existence of the integrating measure will be
completely characterized as an algebraic property of multimoment sequences, Thm. \ref{strong}.
Finally, as an application of the previous results we will provide a new
characterization of a class of weakly harmonizable stochastic processes in section \ref{stochastic}.


\section{The weak multilinear Haussdorff moment problem and a multilinear Riesz theorem }

A polymeasure $\gamma$ on the $\sigma$-algebras
$\Sigma_1, \ldots, \Sigma_n$ is a separately
$\sigma$-additive function on the cartesian product of $\Sigma_1,
\ldots, \Sigma_n$  \cite{Do87} (we will consider here only real or complex polymeasures).
The variation of the polymeasure $\gamma$ is defined as the set
function $v(\gamma)\colon \Sigma_1 \times \cdots \times \Sigma_n
\to [0, +\infty]$:
 $$ v(\gamma ) (A_1, \ldots, A_n) = \sup \left\{ \sum_{k_1}^{r_1} \cdots  \sum_{k_n}^{r_n} | \gamma (A_1^{k_1}, \ldots, A_n^{k_n}) | \right\} ,$$
 where the supremum is taken over all finite partitions $\{ A_l^{k_l} \}_{k_l=1}^{r_l}$
of the set $A_l \in \Sigma_l$.  The semivariation $|| \gamma ||\colon \Sigma_1 \times \cdots
\times \Sigma_n \to [0, +\infty]$ of the polymeasure $\gamma$ is defined as:
\begin{equation}\label{semi}
 || \gamma || (A_1, \ldots, A_n) = \sup \left\{ \Big| \sum_{k_1}^{r_1} \cdots  \sum_{k_n}^{r_n} a_1^{k_1} \cdots a_n^{k_n}\gamma (A_1^{k_1}, \ldots, A_n^{k_n}) \Big| \right\} ,
 \end{equation}
where the supremum is taken over all finite partitions
$\{A_l^{k_l} \}_{k_l=1}^{r_l}$ of the set $A_l \in \Sigma_l$ and
all collections of numbers $\{ a_l^{k_l} \}_{k_l=1}^{r_l}$ such that
$|a_l^{k_l} | \leq 1$. In the linear case $n=1$ the semivariation
and variation of a measure coincide.

An integral denoted as $\int f_1 \otimes \ldots \otimes  f_n \, \gamma  $,
can be constructed for polymeasures of finite semivariation for families of bounded $\Sigma_k$-measurable  scalar
functions $f_k$,  by taking the limits of the integrals of $n$-tuples
of simple functions  uniformly converging to the $f_k$'s \cite{Do87}.  This integral coincides for
compact sets with the integral discussed by Morse and Transue \cite{Mo56} whose main
properties were reviewed in \cite{Ch83}.  Among them we must point it out an extension
of Lebesgue dominated convergence theorem (see example 2.5 and the comments below, Thm. 2.8 and Cor. 2.9 (iii) in \cite{Ch83}).

Let $\mathrm{Bo}(K_l)$ denote the Borel
$\sigma$-algebra on the compact space $K_l$. A polymeasure
$\gamma$ on the product of the $\sigma$-algebras $\mathrm{Bo}(K_1) \times
\cdots \times \mathrm{Bo}(K_n)$ is said to be regular if for any Borel
subsets $A_l \subset K_l$, $l \neq k$, the set function:
$$\gamma_k (A) = \gamma (A_1, \dots, A_{k-1},A,A_{k+1},\dots,A_n)$$
is a signed Radon measure on $K_l$, $l = 1,\ldots, n$.    

The space of
regular coun\-ta\-bly additive polymeasures on $\mathrm{Bo}(K_1)
\times \cdots \times \mathrm{Bo}(K_n)$ will be denoted by
$\mathrm{rcapm}(\mathrm{Bo}(K_1), \ldots, \mathrm{Bo}(K_n))$ which
is a Banach space equipped with the semivariation norm. On the
other hand, we will denote by $\Lie^n(C(K_1), \dots, C(K_n) ; \mathbb{R})$ the space of con\-ti\-nuous scalar
$n$-linear maps on the Banach spaces $C(K_l)$.
Then there exists a natural isometric isomorphism between
$\Lie^n(C(K_1),\dots,C(K_n);\mathbb{R})$, and the space of regular countably additive
polymeasures $\mathrm{rcapm}(\mathrm{Bo}(K_1), \ldots, $ $\mathrm{Bo}(K_n))$
 \cite{Bo98}.

We will use throughout the rest of this paper a consistent
multi--index notation.  We introduce the symbols $\nabla^\br \mu_\bs = \nabla_1^{r_1}\nabla_2^{r_2}\cdots \nabla_n^{r_n} \mu_{s_1\ldots s_n}$, where $\nabla_l$ denotes the difference operator on the $l$th component, $\nabla_l \mu_\bk = \mu_\bk - \mu_{\bk + \mathbf{1}_l}$.  Then we obtain easily:
$$ \nabla^\br \mu_\bs = \sum_{\bl = 0}^\br (-1)^{|\bl |} \binom{\br}{\bl} \mu_{\bs+\bl} $$
where $\br, \bs, \bl$ denote multindexes of length $n$.  We define the Bernstein coefficients $\lambda_{(\bk; \bm)}$ of a $n$--linear functional $L$ as:
\begin{equation}\label{lambda_lambda}
L\left(\binom{k_1}{ m_1}t^{m_1}(1-t)^{k_1-m_1}, \ldots ,\binom{k_n}{m_n}t^{m_n}(1-t)^{k_n-m_n}\right) = \lambda_{(\bk; \bm)} .
\end{equation}
If $L$ is the functional defined by $\mu_\bk$, then $ \lambda_{(\bk; \bm) } = \binom{\bk}{\bm} \nabla^{\bk-\bm} \mu_\bm $.

A function $\mu$ on $\mathbb{N}_0^n$ is called completely monotone if $\nabla^\br \mu \geq 0$ for all $\br$.  Because of eq. (\ref{lambda_lambda}) this is equivalent to the positivity of the functional $L_\mu$, hence to the existence of a Radon measure solving the classical Hausdorff problem of moments and to the $\tau$--positivity of the function $\mu$ (see \cite{Be84} for a thorough discussion of these results, Thm. 4.6.4).  We will introduce now two notions of uniform boundedness for a multi--index sequence $\mu_\bk$ that will characterize the solutions of the classical and weak multilinear Hausdorff problems in the situation where the function $\mu$ is not completely monotone.

\begin{defin}\label{bounded}  A multi--index sequence $\mu_{\bk}$ is said to be bounded
with constant $C>0$ if:
\begin{equation}\label{bound}
 \sum_{\bm = 0}^\bk |\lambda_{(\bk;\bm)}| \leq C, \quad \quad \forall \bk \geq 0 ,
 \end{equation}
and $\mu_{\bk}$ is said to be weakly bounded
with constant $C>0$ if:
\begin{equation}\label{weak_bound}
 \left|\sum_{\bm = 0}^\bk \ba_\bm^\bk \lambda_{(\bk;\bm)}\right| \leq C, \quad \quad \forall \bk\geq 0 ,
 \end{equation}
where $\ba_\bm^\bk = a_{m_1}^{k_1}\cdots  a_{m_n}^{k_n} $, for all $a_{m_l}^{k_l}$ such that $|a_{m_l}^{k_l}| \leq 1$, $l = 1, \ldots, n$.
\end{defin}

It is clear from the definitions that condition (\ref{bound}) implies  (\ref{weak_bound}).

\begin{thm}\label{cont_multi}
A multi--index sequence $\mu_\bk$, $\bk \in \N_0^n$ solves the weak Hausdorff problem if and only if it
is weakly bounded.
\end{thm}

{\parindent 0cm  \bf Proof:}
If $\mu_\bk$ is a solution of the weak Hausdorff multilinear moment problem, eq. (\ref{weak_problem}), then there exists a continuous multilinear functional $L$
on $C[0,1]$ such that $L(\bt^\bk) = \mu_\bk$.    Let $\ba_\bm^\bk = a_{m_1}^{k_1}\cdots  a_{m_n}^{k_n} $, $\bk \geq 0$, with $|a_{m_l}^{k_l}| \leq 1$, $l = 1, \ldots, n$. Then:
\begin{eqnarray}  \left| \sum_{\bm = 0}^\bk \ba_\bm^\bk \lambda_{(\bk;\bm )}  \right| &=&
 \left| \sum_{\bm = 0}^\bk \ba_\bm^\bk L(\lambda_{(k_1,m_1 )}(t), \ldots,  \lambda_{(k_n,m_n )} (t)) \right|  \nonumber \\ & = & \nonumber
 \left|  L\left(\sum_{m_1 = 0}^{k_1} a_{k_1}^{m_1}\lambda_{(k_1,m_1 )}(t), \ldots, \sum_{m_n = 0}^{k_n}a_{k_n}^{m_n} \lambda_{(k_n,m_n )}(t) \right )\right|  \\ & \leq & || L || \prod_{l = 1}^n || \sum_{m_l = 0}^{k_l} a_{k_l}^{m_l}\lambda_{(k_l,m_l )}(t) ||_\infty  \nonumber \end{eqnarray}
Moreover, because $\left| \sum_{m_l = 0}^{k_l} a_{k_l}^{m_l}\lambda_{(k_l,m_l )}(t) \right|  \leq   \sum_{m_l = 0}^{k_l}\binom{k_l}{m_l} t^{k_l -m_l}(1-t)^{m_l} = 1$, we reach the conclussion.

\medskip

Conversely, if we assume that the multi--index sequence $\mu_\bk$ of length is weakly bounded,  we will show by
induction on $n$ that $ | L (p_1, \ldots, p_n ) | \leq 2^n C || p_1||_\infty\cdots || p_n||_\infty$, for any family of polynomials $p_1, \ldots, p_n$.

For $n=1$, because of \eqref{linear_bound} the weakly bounded condition for sequences is equivalent to the condition of boundedness.

We will assume that if $L'$ is a $(n-1)$--multilinear
functional associated to the bounded $(n-1)$--multi--index
sequence $\mu'_{\bk'}$ with bounding constant $C'$, then
$$| L'(q_1, \ldots, q_{n-1} ) | \leq 2^{n-1} C' || q_1||_\infty \cdots || q_{n-1} ||_\infty ,$$
for any family of polynomials $q_1, \ldots, q_{n-1}$.
Let  $L$ be the multilinear functional
associated to the weakly bounded $n$--multi--index sequence $\mu_{\bk}$.
We denote by $L_p$ the $(n-1)$--multilinear functional obtained by
fixing the $n$th argument of $L$ to be the polynomial $p$, i.e., $L_p(q_1, \ldots , q_{n-1}) = L (q_1, \ldots, q_{n-1}, p)$.

Notice that if $p(t)$ is a polynomial of degree $r$ then \cite{Sh70}:
$$B_N(p)(t) = p(t) + S_N(t) = p(t) +  \sum_{l=1}^{r-1}\frac{p_{r,l}(t)}{N^l} ,$$
where $B_N(p)$ denotes the $N$th Bernstein polynomial of the function $p(t)$ and  $p_{r,l}$ are polynomials of degree less than or equal to $r$, not depending on $N$.
Denoting by $a_{\bm'}^{\bk'} = a_{m_1}^{k_1}\cdots a_{m_{n-1}}^{k_{n-1}}$ with $|a_{m_l}^{k_l}| \leq 1$,  it is clear that:
 \begin{eqnarray*}  &&
 \left| \sum_{\bm' = 0}^{\bk'} a_{\bm'}^{\bk'} L_{B_N(p)} (\lambda_{k_1,m_1}(t_1), \ldots, \lambda_{k_{n-1},m_{n-1}}(t_{n-1})) \right|  \leq \\
&& \leq  \left| \sum_{\bm' = 0}^{\bk'} \sum_{m_n= 0}^N a_{\bm'}^{\bk'} p( m_n/N) L  (\lambda_{k_1,m_1}(t_1), \ldots, \lambda_{k_{n-1},m_{n-1}}(t_{n-1}), \lambda_{N,m_n}(t_n)) \right| \leq \\
&& \leq || p ||_\infty \left| \sum_{\bm = 0}^{\bk}  a_{\bm}^{\bk} L( \lambda_{k_1,m_1}(t_1), \ldots, \lambda_{k_{n-1},m_{n-1}}(t_{n-1}), \lambda_{k_n,m_n}(t_n)) \right|
\end{eqnarray*}
with $\bk = (k_1, \ldots, k_{n-1},N)$,  $\bm = (m_1, \ldots, m_{n_1},m_n)$ and $a_{\bm}^{\bk} = a_{\bm'}^{\bk'} ( p( m_n/N) / ||p||_\infty ).$
Hence,
\begin{eqnarray*}
 \left| \sum_{\bm' = 0}^{\bk'} a_{\bm'}^{\bk'} L_{B_N(p)} (\lambda_{k_1,m_1}(t_1), \ldots, \lambda_{k_{n-1},m_{n-1}}(t_{n-1})) \right| & \leq & \\ \leq  || p ||_\infty \left| \sum_{\bm = 0}^{\bk} a_{\bm}^{\bk}   \lambda_{(\bk;\bm)} \right| & \leq &  C || p ||_\infty
 \end{eqnarray*}
because $\mu_{\bk}$ is weakly bounded with constant $C$.  Then the $(n-1)$--multimoment sequence defined by $L_{B_N(p)}$ is weakly bounded with bound $C||p||_\infty$, and by the induction hypothesis we obtain:
$$ | L_{B_N(p)} (p_1,\ldots, p_{n-1} ) | \leq 2^{n-1} C || p_1||_\infty \cdots || p_{n-1} ||_\infty ||p||_\infty . $$

Similarly, we consider now $S_N(t) =  \sum_{l=1}^{r-1}\frac{p_{r,l}(t)}{N^l} $.
If the polynomials $p_{r,l}$ have the form, $p_{r,l}(t) = \sum_{j = 0}^{r} a_{lj} t^j$, then by choosing:
\begin{equation}\label{conditions2}
a = \max \{ |a_{lj} |  \}, \quad N \geq \frac{a(r-1) (r + 1)}{|| p ||_\infty} ,
\end{equation}
and using the notations above, we will get:
\begin{eqnarray} \label{last}
& &\left| \sum_{\bm' = 0}^{\bk'} a_{\bm'}^{\bk'} L_{S_N(p)}(\lambda_{k_1,m_1}(t_1), \ldots, \lambda_{k_{n-1},m_{n-1}}(t_{n-1})) \right| \leq \\ \nonumber
& & \leq \left| \sum_{\bm' = 0}^{\bk'} a_{\bm'}^{\bk'} \sum_{l=1}^{r-1} L(\lambda_{k_1,m_1}(t_1), \ldots, \lambda_{k_{n-1},m_{n-1}}(t_{n-1}), \frac{p_{r,l}(t_n)}{N^l}) \right| = \\  \nonumber
& & =  \left| \sum_{\bm' = 0}^{\bk'} a_{\bm'}^{\bk'} \sum_{j= 0}^r \sum_{l=1}^{r-1} \frac{a_{lj}}{N^l} L(\lambda_{k_1,m_1}(t_1), \ldots, \lambda_{k_{n-1},m_{n-1}}(t_{n-1}), t_n^j) \right| \leq \\ \nonumber
& & \leq  \frac{|| p ||_\infty}{r+1}  \sum_{j= 0}^r  \left| \sum_{\bm' = 0}^{\bk'} a_{\bm'}^{\bk'}\left( \frac{r+1}{||p||_\infty}\sum_{l=1}^{r-1} \frac{a_{lj}}{ N^l} \right)L(\lambda_{k_1,m_1}(t_1), \ldots \right. \\  \nonumber
&  & \quad \left. \ldots , \lambda_{k_{n-1},m_{n-1}}(t_{n-1}), \lambda_{j,j}(t_n)) \right| .
\end{eqnarray}
If we denote by $a_j^j$ the quantity $\frac{r+1}{||p||_\infty}\sum_{l=1}^{r-1} \frac{a_{lj}}{ N^l}$,  conditions eq. (\ref{conditions2}),
imply that $| a_j^j| \leq 1$.  We will consider now the numbers $a_{m_n}^j = 0$ for all $0 \leq m_n < j$. With these definitions the last term in the sequence of inequalities (\ref{last}), can be written as:
\begin{eqnarray*}
& & = \frac{|| p ||_\infty}{r+1}  \sum_{j= 0}^r  \left| \sum_{\bm' = 0}^{\bk'}  \sum_{m_n= 0}^j a_{\bm'}^{\bk'}a_{m_n}^j L(\lambda_{k_1,m_1}(t_1), \ldots \right. \\
& & \hskip 5cm  \left. \ldots , \lambda_{k_{n-1},m_{n-1}}(t_{n-1}), \lambda_{j,m_n}(t_n)) \right| = \\ \nonumber
& & = \frac{|| p ||_\infty}{r+1}  \sum_{j= 0}^r  \left| \sum_{\bm = 0}^{\bk}  a_{\bm}^{\bk } L(\lambda_{k_1,m_1}(t_1), \ldots, \lambda_{k_{n-1},m_{n-1}}(t_{n-1}), \lambda_{j,m_n}(t_n)) \right|
\end{eqnarray*}
with $\bk = (k_1, \ldots, k_{n-1},j)$ and $\bm = (m_1, \ldots, m_n)$.  Hence finally we obtain:
\begin{eqnarray*}
& & \left| \sum_{\bm' = 0}^{\bk'} a_{\bm'}^{\bk'} L_{S_N(p)}(\lambda_{k_1,m_1}(t_1), \ldots, \lambda_{k_{n-1},m_{n-1}}(t_{n-1})) \right| \leq \\
& & \leq  \frac{|| p ||_\infty}{r+1}  \sum_{j= 0}^r  \left| \sum_{\bm = 0}^{\bk}  a_{\bm}^{\bk } \lambda_{(\bk ; \bm )}\right| \leq  \frac{|| p ||_\infty}{r+1}  \sum_{j= 0}^r  C = C || p ||_\infty
\end{eqnarray*}
and the sequence of multimoments  $\mu''_{\bk''} = L_{S_N(p)}(\lambda_{k_1,m_1}(t_1), \ldots, \lambda_{k_{n-1},m_{n-1}})$
is  weakly bounded with constant $C||p||_\infty$.

We conclude the argument by using the induction hypothesis and computing:
\begin{eqnarray*}  | L (p_1, \ldots, p_n ) | &=& |L_{p_n} (p_1, \ldots, p_{n-1})| \\
&\leq &   |L_{B_N(p_n)} (p_1, \ldots, p_{n-1})|  +  |L_{S_N(p_n)} (p_1, \ldots, p_{n-1})|  \\
&\leq& 2^n C||p_1||_\infty \cdots || p_n ||_\infty \quad\quad \quad\mbox{\hskip 3cm}\Box
\end{eqnarray*}

\bigskip

It is clear that if the multi--index sequence $\mu_{\bk}$ is nonnegative, i.e., $\mu_{\bk}$ are positive or zero real numbers for all $\bk$, then
the sequence $\mu_{\bk}$ is weakly bounded iff is bounded because,
$$\sum_{\bm = 0}^\bk |\lambda_{(\bk;\bm)}| = \sum_{\bm = 0}^\bk \lambda_{(\bk;\bm)} \leq \sup_{\ba_\bm^\bk, |a_{m_l}^{k_l}| \leq 1} \left|\sum_{\bm = 0}^\bk \lambda_{(\bk;\bm)} \right|\leq C .$$
Moreover under these circumstances, it is simple to see that the total variation of the polymeasure $\gamma$ determined by $\mu_{\bk}$ is finite, hence 
the polymeasure $\gamma$ determines a Radon measure on $[0,1]^n$ \cite{Bo01}.
Thus we have obtained a particular instance of the fact that Radon bimeasures on Hausdorff spaces are extensions of Radon measures on the product spaces (see \cite{Be84}, Thm. 1.1.10).

\section{The strong Hausdorff multilinear moment problem: orthogonally additive polynomials}

Given a multi--index sequence $\mu_{\bk}$, $\bk =
(k_1, \cdots, k_n)$, $k_l = 0,1,\ldots$, $l = 1, \ldots, n$, we will call it H\"ankel if 
$\mu_{\bk  + \mathbf{1}_l} = \mu_{\bk  + \mathbf{1}_{l+1}}$, where the multi-index $\mathbf{1}_l $ is defined as $(\mathbf{1}_l)_j  = \delta_{lj}$,  for all $l = 1, \ldots, n$.

\begin{thm}\label{strong}  Let  $\mu_{\bk}$ be a multi--index sequence solving the classical Hausdorff problem of moments.  Then $\mu_{\bk}$ solves the strong Hausdorff problem of moments if and only if  $\mu_{\bk}$ is H{\"a}nkel.
\end{thm}

{\parindent 0cm \bf Proof:}   Consider the $n$--linear
functional defined on the space $\mathcal{P}$ of real polynomials
on $I = [0,1]$ by the multi--index sequence $\mu_{\bk}$.
Because the multi--index sequence $\mu_{\bk}$ is bounded, then
$L$ can be extended to $C(I)$ (Thm. \ref{cont_multi}).    We shall
denote such extension with the same symbol $L$. 

The homogeneous polynomial $P_L$ determined by $L$ is
orthogonally additive.  To prove it we notice that $L(f_1, \ldots,
g\cdot f_l , f_{l+1}, \ldots, f_n) = L (f_1, \ldots, f_l, g\cdot
f_{l+1}, \ldots, f_n)$ for all $f_1, \ldots, f_n, g\in C(I)$.  In
fact we can construct a sequence of polynomials $p_{m_l}$, $q_m$
converging uniformly  to $f_l$ and $g$ respectively ($l = 1,\ldots,
n$) on $I$, hence because $\mu_{\bk}$
is H{\"a}nkel, we have:
$$
 L(p_{m_1}, \ldots, q_m\cdot p_{m_l} , p_{m_{l+1}}, \ldots, p_{m_n}) =
 L (p_{m_1}, \ldots, p_{m_l} , q_m\cdot p_{m_{l+1}}, \ldots, p_{m_n}) ,
 $$
and the conclusion follows because of the continuity of $L$.

Now suppose we have two disjoint positive functions $f,g$ on $C(I)$, $|f|\wedge |g| = 0$.  We compute:
\begin{eqnarray*}
P_L(f+g) &=& L(f+g, \ldots, f+g) = \sum_{r \geq 0}  \binom{n}{r}  L(f, {}_{\cdots}^{n-r},f,g, {}_{\cdots}^r ,g ) = \\
&=& L(f,\ldots, f) + \sum_{r =1}^{n-1}  \binom{n}{r}  L (1,f,{}_{\cdots}^{n-r-1},f,g,{}_{\cdots}^{r-1},g,f\cdot g) \\
&+& L(g,\ldots, g) = P_L(f) + P_L(g),
\end{eqnarray*}
because $f\cdot g = 0$.

Using the representation theorem for orthogonally additive polynomials on Banach lattices \cite{Be06} and because the $n$-concavification of the Banach lattice $E = C(I)$
coincides  with itself, this is $C(I)_{(n)} = C(I)$ (see also \cite{Ca06} and \cite{Pe05}), the
polynomial $P_L$ defines a bounded linear functional
$T\colon C(I) \to \R$,
$$T(f^n) = P_L(f) = L(f,\cdots, f)$$
and then, by Riesz theorem, there will exists a Radon measure
$\mu$ such that $T(f^n) = \int_I f(t)^n\, \mu(dt)$.   Hence
$L(t^{k_1}, \ldots, t^{k_n}) = \int_I t^{k_1+\cdots+k_n} \, \mu(dt)$, and $ L (f_1, \ldots, f_n ) = \int_I f_1(t)\cdots f_n(t) \, \mu(dt)$. \hfill
$\Box$

\section{The weak bilinear Hausdorff problem of moments and weakly harmonizable stochastic processes}\label{stochastic}

We will use now the characterization of polymeasures with compact support
obtained before to provide a description of a class of weakly harmonizable processes.   
Let us consider an stochastic process $X_t$ modelled on a probability space $(\Omega, \Sigma, P)$ where $\Sigma$ is a $\sigma$--algebra on the set $\Omega$,
$P$ is a probability measure on $\Omega$, and the map $X\colon \mathbb{R} \to L^2(\Omega, P)$, $X_t := X(t)$, is strongly continuous.   We denote by $C(t,t') = E(\bar{X}_t X_t')$ the covariance function of $X_t$.    If the second order process $X_t$ is weakly stationary, i.e., there exists an univariate continuous function $\Phi$ such that $C(t,t') = \Phi(t'-t)$, Cram\'er--Kolmogorov's theorem shows that there exists a stochastic measure $\xi$ on $\mathbb{R}$ with values on $L^2(\Omega,P)$ such that the process $X_t$ is the Fourier transform of $\xi$:
\begin{equation}\label{harmonic}
X_t = \int_{\mathbb{R}} e^{its} \, \xi (ds) ,
\end{equation}
Moreover if $A,B$ are two Borel sets on $\mathbb{R}$, then
\begin{equation} \label{orthogonality}
\langle \xi (A), \xi(B) \rangle_{L^2(\Omega, P)} = \mu(A\cap B) .
\end{equation}
Such processes are called (strongly) harmonizable.   An important class of second order stochastic processes $X_t$ that admit generalized harmonic representations are the so called weakly harmonizable and they satisfy:
$$ C(t,t') =  \int_{\mathbb{R}\times \mathbb{R}} e^{-its} \otimes e^{it's'} \, \gamma (ds, ds') ,$$
where $\gamma$ is a positive definite bimeasure, this is: 
\begin{equation}\label{positivity}
\gamma (A,B) = \overline{\gamma (B,A)},  \quad   \sum_{i,j=1}^r \bar{a}_i a_j \gamma (A_i, A_j) \geq 0\, ,
\end{equation}
for all families of complex numbers $a_i$ and Borel sets $A,B, A_i$ on $\mathbb{R}$.  
Then if $X_t$ is a weakly harmonizable process then there exists a harmonic representation of the form eq. (\ref{harmonic}) for them, where now the orthogonality condition (\ref{orthogonality}) is replaced by $\langle \xi (A), \xi(B) \rangle_{L^2(\Omega, P)} = \gamma(A , B)$
 for any Borel sets $A,B$ (see for instance the review \cite{Ra82}, Thm. 3.2.)

We will consider a complex regular bimeasure $\gamma$ of finite semivariation with support in $[0,1]\times [0,1]$.  Let us call such bimeasures Hausdorff.   We will consider the Fourier--Stieltjes transform of the  bimeasure $\gamma$:
$$ \hat{\gamma} (t,t') = \int_{\mathbb{R}\times \mathbb{R}} e^{-its}\otimes e^{it's'} \, \gamma (ds, ds') .$$
The function $\hat{\gamma}$ is bounded by $\hat{\gamma} (0,0) = \gamma([0,1],[0,1]) = \mu_{00} \leq ||\gamma || <  \infty$ and the extension of Lebesgue's dominated convergence theorem mentioned in section 2.1 shows that the function $\hat{\gamma}$ is analytic in the real plane $(t,t')$ with power series expansion given by:
\begin{equation} \label{fourier_transform}
\hat{\gamma} (t,t')  = \sum_{n,m\geq 0} (-1)^n\frac{i^{n+m}}{n!m!} \mu_{nm} t^n t'^m ,
\end{equation}
where the coefficients $\mu_{nm}$, $n,m\geq 0$ are the moments of the bimeasure $\gamma$.   
Moreover the sequence of moments $\mu_{nm}$ is weakly bounded because of Thm. \ref{cont_multi}.
Now a simple argument shows that these conditions characterize completely the Fourier-Stieltjes transform of Hausdorff polymeasures.

A weakly harmonizable second order stochastic process $X_t$ such that the support of the stochastic measure $\xi$ defining it is contained in the interval $[0,1]$ will be called Hausdorff.   Notice that in such a case because of eq. (\ref{harmonic}), the support of the corresponding bimeasure $\gamma$ will be contained in $[0,1]\times [0,1]$. 
Now if we are given a arbitrary collection of complex numbers $a_1, \ldots, a_r$, and we compute $\sum_{l,k = 1}^r \bar{a}_l a_k \hat{\gamma}(t_l,t_k)$  for a positive definite bimeasure $\gamma$, we obtain:
$$\sum_{l,k = 1}^r \bar{a}_l a_k \hat{\gamma}(t_l,t_k) = \int_{\mathbb{R}\times \mathbb{R}} \sum_{l=1}^r \bar{a}_le^{-it_ls} \otimes \sum_{k = 1}^ra_ke^{it_ks'} \, \gamma (ds, ds') \geq 0.$$
and we conclude that  the function $\hat{\gamma}$ is a positive definite kernel.  Hence the class of analytic positive definite kernels described above are just the covariance functions of Hausdorff weakly harmonizable stochastic processes.

\begin{cor}  A function of two real variables $\Phi (t,s)$ is the covariance function of a second order weakly harmonizable Hausdorff process $X_t$ if and only if is an  analytic positive definite kernel on $\mathbb{R}^2$ such that the multi--index sequence $\mu_{nm} = \partial^{n+m}\Phi/\partial t^n \partial s^m (0,0)$, $n,m \in \mathbb{N}_0$ is weakly bounded.  Moreover, the stochastic process $X_t$ will be weakly stationary if and only if the multimoment sequence $\mu_{nm}$ is H\"ankel, i.e., $\mu_{n+1,m} = \mu_{n,m+1}$ for all $n,m$. 
\end{cor}

\begin{rem}  Notice that the positivity condition can be dispensed with as it follows from the previous discussion that an  analytic function $\Phi(s,t)$ wiill have the form (\ref{fourier_transform}) for a Hausdorff bimeasure $\gamma$ iff its sequence of moments is weakly bounded.  However, unless the positivity condition (\ref{positivity}) is satisfied it is not possible to reconstruct a Hilbert space where the stochastic process would be represented. 
\end{rem}



\subsection*{Acknowledgements}
AI and JGL have been partially supported by Spain MICIN Project MTM2010-21186-C02-2 and QUITEMAD project.
The second author was partially supported by FPU-MEC AP-2004-4843
Grant and by the ``Programa de formaci{\'o}n del profesorado
universitario del MEC" and by Project MTM 2006-03531.  The authors would like to thank F. Bombal and I. Villanueva for their
technical support with polymeasures.

\end{document}